\documentclass[11pt]{amsart}
\usepackage[active]{srcltx}
\usepackage{calc,amssymb,amsthm,amsmath,amscd, eucal,ulem}
\usepackage{alltt}
\usepackage{fullpage}
\RequirePackage[dvipsnames,usenames]{color}

\input{kmacros3.sty}
\input{mabliautoref.sty}
\input{xy}
\xyoption{all}
\usepackage{tikz}

\numberwithin{equation}{theorem}
\def\ff{{\bf f}}
\def\XX{{\bf X}}
\def\YY{{\bf Y}}

\DeclareMathOperator{\edim}{edim}

\DeclareMathOperator{\homgp}{Hom}

\DeclareMathOperator{\sheafhom}{\scr{H}{\kern -2pt om}}

\DeclareMathOperator{\ffield}{frac} 

\usepackage{fullpage}

\usepackage{setspace}
\usepackage{hyperref}

\usepackage{enumerate}

\usepackage{graphicx}

\usepackage[all,cmtip]{xy}
\usepackage{verbatim}

\theoremstyle{theorem}

\begin{document}
\title{RationalMaps, a package for Macaulay2}
\author{C.J. Bott}
\author{S. Hamid Hassanzadeh}
\author{Karl Schwede}
\author{Daniel Smolkin}
\address{Department of Mathematics\\Mailstop 3368\\Texas A\&M University\\College Station, TX 77843-3368}
\email{cbott2@math.tamu.edu}
\address{Department of Mathematics\\Federal University of Rio de Janeiro\\Brazil}
\email{hamid@im.ufrj.br}
\address{Department of Mathematics\\ University of Utah\\ Salt Lake City\\ UT 84112}
\email{schwede@math.utah.edu}
\address{Department of Mathematics\\ University of Michigan\\ Ann Arbor\\ UT 48109}
\email{smolkind@umich.edu}

\thanks{The second named author was supported by CNPq-bolsa de Produtividade and by the MathAmSud project ``ALGEO''}
\thanks{The third named author was supported in part by the NSF FRG Grant DMS \#1265261/1501115, NSF CAREER Grant DMS \#1252860/1501102, NSF Grants DMS \#1840190 and DMS \#2101800.}
\thanks{The fourth named author was supported in part by the NSF FRG Grant DMS \#1265261/1501115, NSF CAREER Grant DMS \#1252860/1501102 and NSF Grant DMS \#1801849.}

\begin{abstract}
This paper describes the {\tt RationalMaps} package for Macaulay2.  This package provides functionality for computing several aspects of rational maps.
\end{abstract}
\maketitle

\section{Introduction}

This package aims to compute several things about rational maps between varieties.  In particular, this package will compute
\begin{itemize}
\item{} The base locus of a rational map.
\item{} Whether a rational map is birational.
\item{} The inverse of a birational map.
\item{} Whether a map is a closed embedding.
\item{} And more!
\end{itemize}
Our functions have numerous options which allow them to run much more quickly in certain examples if configured correctly.  Setting the option {\tt Verbosity} to a value $\geq 1$ will mean that functions will provide hints as to the best ways to run them.  This paper discusses {\tt RationalMaps} version 1.0.

A rational map $\mathfrak{F}:X\subseteq \P^n\dasharrow Y \subseteq \P^m$ between projective varieties is presented by $m+1$ forms $\ff=\{f_0,\ldots f_m\}$ of the same degree in the coordinate ring of $X$, denoted by  $R$.
The idea of looking at the syzygies of the forms $\ff$ to detect the geometric properties of  $\mathfrak{F}$ goes
back at least to \cite{HulekKatzSchreyer} in the case
where $X=\P^n$, $Y = \P^m$ and $m=n$ (see also \cite{SempleTyrrell}). In \cite{RussoSimisCompositio}  this method was developed by Russo and Simis
to handle the case $X=\P^n$ and $m\geq n$. Simis pushed the method further to the study
of general rational maps between two integral projective schemes in
arbitrary characteristic by an extended ideal-theoretic method
emphasizing the role of the Rees algebra associated to the ideal
generated by $\ff$ \cite{SimisCremona}.  Recently,  Doria, Hassanzadeh, and Simis applied these
 Rees algebra techniques to study the birationality of  $\mathfrak{F}$ \cite{DoriaHassanzadehSimisBirationality}.  Our core functions, particularly those related to computing inverse maps, rely heavily on this work.
\vskip 12pt
\noindent
\emph{Acknowledgements:}  Much of the work  on this package was completed during the Macaulay2 workshop held at the University of Utah in May 2016.  We especially thank David Eisenbud, Aron Simis, Greg Smith, Giovanni Staglian\`o,  and Mike Stillman for valuable conversations.  We also thank the referees for numerous helpful comments on both the package and the paper.

\section{Base Loci}

We begin with the problem of computing the base locus of a map to projective space. Let $X$ be a projective variety over any field $k$ and let $\mathfrak{F}: X \to \P_k^m$ be a rational map from $X$ to projective space. Then we  choose a representative $(f_0, \cdots, f_m)$ of $\mathfrak{F}$, where each $f_i$ is the $i^{\textrm{th}}$ coordinate of $\mathfrak{F}$. A priori, each $f_i$ is in $K = \ffield R$, where $R$ is the coordinate ring of $X$. However, we can get another representative of $\mathfrak{F}$ by clearing denominators. (Note that this does not enlarge the base locus of $\mathfrak{F}$ since $\mathfrak{F}$ is undefined whenever the denominator of any of the $f_i$ vanishes.) Thus we  assume that $f_i\in R$ for all $i$, and that all the $f_i$ are homogeneous of the same degree.

In this setting, one might naively think that the map $\mathfrak{F}$ is undefined exactly when all of the $f_i$ vanish, and thus the base locus is the vanishing set of the ideal $(f_0, \cdots, f_m)$. However, this  yields a base locus that's too big.  Indeed, to find the base locus of a rational map, we must consider all possible representatives of the map and find where none of them are defined. To do this, we use the following result.


%

\begin{proposition}\textnormal{\cite[Proposition 1.1]{SimisCremona}}
  Let $\mathfrak{F}: X \dashrightarrow \P^m$ be a rational map and let $\textbf{f} = \left\{ f_0, \dots, f_m \right\}$ be a representative of $\mathfrak{F}$ with $f_i\in R$ homogeneous of degree $d$ for all $i$. Set $I  = (f_0, \cdots, f_m)$. Then the set of such representatives of $\mathfrak{F}$ corresponds bijectively to the homogeneous vectors in the rank 1 graded $R$-module $\homgp_R(I, R) \cong (R :_K I)$.
  \label{lemma:repsOfRatMap}
\end{proposition}

The bijection comes from multiplying our fixed representative $\textbf{f}$ of $\mathfrak{F}$ by $h \in (R :_K I)$.
Now, in the setting of \autoref{lemma:repsOfRatMap}, let
  \[
    \bigoplus_s R(-d_s) \xrightarrow{\varphi} R(-d)^{m+1} \xrightarrow{[f_0, \cdots, f_m]} I \to 0
  \]
be a free resolution of $I$. Then we get
\[
  0 \to \homgp_R(I, R) \to \left( R(-d)^{m+1} \right)^\vee \xrightarrow{\varphi^t} \left( \bigoplus_s R(d_s) \right)^\vee
\]
where $\varphi^t$ is the transpose of $\varphi$ and $R^\vee$ is the dual module of $R$. Thus, we get that $\homgp_R(I,R) \cong \ker \varphi^t$, and so each representative of $\mathfrak{F}$ corresponds to a vector in $\ker \varphi^t$. The correspondence takes a representative $(hf_0, \cdots, hf_m)$ to the map that multiplies vectors in $R^{m+1}$ by $[hf_0, \cdots, hf_m]$ on the left.

The base locus of $\mathfrak{F}$ is the intersection of the sets $V(f^i_0, \cdots, f^i_m)$ as $\mathbf{f}^i = (f^i_0, \cdots, f^i_m)$ ranges over all the representatives of $\mathfrak{F}$. The above implies that this is the same as the intersection of the sets $V(w^i_0,\cdots, w^i_m)$ as $\mathbf{w}^i = (w^i_0, \cdots, w^i_m)$ ranges over the vectors in $\ker \varphi^t$. Now, given any $a, f, g\in R$, we have $V(af) \supseteq V(f)$ and $V(f + g) \supseteq V(f)\cap V(g)$. Thus, it's enough to take a generating set $\mathbf w^1, \cdots, \mathbf w^n$ of $\ker \varphi^t$  and take the intersection over this generating set.

The base locus of $\mathfrak{F}$ is then the variety cut out by the ideal generated by all the entries of all of the $\mathbf w^i$. Our function {\tt baseLocusOfMap}  returns this ideal.  It can be applied either to our new Type {\tt RationalMapping} or to a {\tt RingMap} between the homogeneous coordinate rings which represents the rational map.
{\scriptsize\color{blue}
\begin{verbatim}
    i1 : loadPackage "RationalMaps";

    i2 : R = QQ[x,y,z];    

    i3 : f = rationalMapping(R, R, {x^2*y, x^2*z, x*y*z})    

                                    2    2
    o3 = Proj R - - - > Proj R   {x y, x z, x*y*z}    

    o3 : RationalMapping    

    i4 : baseLocusOfMap(f)    

    o4 = ideal (y*z, x*z, x*y)    

    o4 : Ideal of R    
\end{verbatim}
}
{\normalsize}
If the \verb=SaturateOutput= option is set  to {\tt false}, our function {\tt baseLocusOfMap} will not saturate the output.

\section{Birationality and Inverse Maps}

Again, a rational map $\mathfrak{F}:X\subseteq \P^n\dasharrow Y \subseteq \P^m$ between projective spaces is defined by $m+1$ forms $\ff=\{f_0,\ldots f_m\}$ of the same degree in the coordinate ring of $X$, denoted by  $R$. $R$ is a standard graded ring in $n+1$ variables. Here we assumed the varieties are defined over a field $k$ and  $\dim R\geq 1$.
Our goal is to find a ring theoretic criterion for birationality and, on top of that, to find the inverse of a rational map.  To do this, we study the Rees algebra of the ideal $I=(\ff)$ in $R$.
To that end set $R\simeq k[x_0,\ldots, x_n]=k[\XX]/\mathfrak{a}$ with $k[\XX]=k[X_0,\ldots, X_n]$ and $\mathfrak{a}$ a homogeneous ideal. The Rees algebra is defined by the polynomial relations among $\{f_0,\ldots f_m\}$ in $R$. To this end, we consider the polynomial extension $R[\YY]=R[Y_0,\ldots,Y_m]$. To keep track of the variables by degrees, we set the standard  bigrading $\deg(X_i)=(1,0)$ and $\deg(Y_j)=(0,1)$. Mapping $Y_j\mapsto f_jt$ yields
a presentation $R[\YY]/\mathcal{J}\simeq {\mathcal R}_R((\ff))$, with $\mathcal{J}$ a bihomogeneous {\em presentation
ideal}.
$\mathcal{J}$ is a bigraded ideal that depends only on the rational map defined by $\ff$
and not on this particular representative.


$${\mathcal J}=\bigoplus_{(p,q)\in \mathbb{N}^2} {\mathcal J}_{(p,q)},$$
where ${\mathcal J}_{(p,q)}$ denotes the $k$-vector space of forms of bidegree $(p,q)$.
Every piece of this ideal contains information about the rational map. For example ${\mathcal J}_{0,*}$ determines the dimension of the image of the map.
For birationality, the following bihomogeneous piece is  important:
\[
    {\mathcal J}_{1,*}:=\bigoplus_{q\in\mathbb{N}} {\mathcal J}_{1,q}
\]
with ${\mathcal J}_{1,q}$ denoting the bigraded piece of ${\mathcal J}$ spanned by the forms of bidegree
 $(1,q)$ for all $q\geq 0$. Now, a form of bidegree $(1,*)$ can be written as $\sum_{i=0}^n Q_i(\YY)\,x_i$, for suitable homogeneous $Q_i(\YY)\in  R[\YY]$
of the same degree.

One then goes to construct a matrix that  measures the birationality of the map. The first step is to lift the polynomials $Q_i(\YY)\in  R[\YY]$ into $k[\XX,\YY]$.
Since the $\{y_0,\cdots,y_m\}$ are indeterminates over $R$, each pair of such representations of the same form gives a syzygy of $\{x_0,\ldots,x_n\}$
with coefficients in $k$.  This is where one must take into attention whether $X\subseteq \P^n$ is minimally embedded  or not. To measure this one can easily check the vector space dimension of ${\mathfrak a}_1$, the degree-1 part of $\mathfrak a$; if it is zero  then $X\subseteq \P^n$ is non-degenerated.



Next, one can pick a minimal set of generators of the ideal $({\mathcal J}_{1,*})$ consisting of a finite number
of forms of bidegree $(1,q)$, for various $q$'s.
Let's assume  $X\subseteq \P^n$ is non-degenerated. Let $\{P_1,\ldots,P_s\}\subset k[\XX,\YY]$ denote liftings of these biforms,
 consider the Jacobian matrix of the polynomials $\{P_1,\ldots,P_s\}$ with respect to $\{x_0,\cdots,x_n\}$. This is a matrix with entries in $k[\YY]$.
Write $\psi$ for the corresponding matrix over $S=k[\YY]/{\mathfrak b}$, the coordinate ring of $Y$. This matrix is called   the \emph{weak Jacobian dual matrix}  associated to
the given set of generators of $({\mathcal J}_{1,*})$.
Note that a weak Jacobian matrix $\psi$ is not uniquely defined due to the lack of uniqueness in the expression of
an individual form and to the choice of bihomogeneous generators. However, it is shown in \cite[Lemma 2.13]{DoriaHassanzadehSimisBirationality} that if the weak Jacobian matrix associated  to one set of bihomogeneous minimal generators of
$({\mathcal J}_{1,*})$ has  rank over $S$ then the weak Jacobian matrix associated to any other
set of bihomogeneous minimal generators of
$({\mathcal J}_{1,*})$ has  rank over $S$ and the two ranks coincide.

\vspace{0.2in}
The  following  criterion is   \cite[Theorem 2.18 ]{DoriaHassanzadehSimisBirationality}. In the package, we consider only the cases where $\XX$ is irreducible i.e. $R$ is a domain.
\begin{theorem} \label{T birationality}Let $X\subseteq \P^n$ be non-degenerate. Then $\mathfrak{F}$ is birational onto $\YY$  if and only if ${\rm rank}(\psi)=\edim(R)-1(=n).$
Moreover
\begin{enumerate}
  \item[{\rm (i)}] We get a representative for the inverse of $\mathfrak{F}$ by taking the coordinates of any homogeneous vector of positive degree
in the {\rm (}rank one{\rm )} null space of $\psi$ over $S$ for which these coordinates generate an ideal containing
a regular element.

\item[{\rm (ii)}] If, further, $R$ is a domain, the representative of $\mathfrak{F}$ in {\rm (i)}
can be taken to be the set of the {\rm (}ordered, signed{\rm )} $(\edim(R)-1)$-minors
of an arbitrary $(\edim(R)-1)\times \edim(R)$ submatrix of $\psi$ of rank $\edim(R)-1$.
\end{enumerate}

\end{theorem}

As expected, the most expensive part of applying this theorem is computing the Rees ideal ${\mathcal J}$. In the package {\tt RationalMaps}, we use {\tt ReesStrategy} to compute the Rees equations. The algorithm is the standard elimination technique. However, we do not use the {\tt ReesAlgebra} package, since verifying birationality according to \autoref{T birationality} only requires  computing a small part of the Rees ideal, namely elements of first-degree $1$. This idea is applied in the {\tt SimisStrategy}. More precisely, if the given map $\mathfrak{F}$ is birational, then the Jacobian dual rank will attain its maximum value of $\edim(R)-1$ after computing the Rees equations up to  degree $(1,N)$ for $N$ sufficiently large. This allows us to compute the inverse map.  The downside of {\tt SimisStrategy} is that if $\mathfrak{F}$ is not birational,  the desired number $N$ cannot be found and the process never terminates. To provide a definitive answer for birationality,  we use {\tt HybridStrategy}, which is a hybrid of {\tt ReesStrategy} and {\tt SimisStrategy}.  The default strategy is {\tt HybridStrategy}.

{\tt HybridLimit} is an option to switch   {\tt SimisStrategy} to  {\tt{ReesStrategy}}, if the computations up to degree $(1, \tt{HybridLimit})$ do not lead to   ${\rm rank}(\psi)=\edim(R)-1$.
The default value for {\tt HybridLimit} is $15$. The change from  {\tt SimisStrategy} to  {\tt ReesStrategy} is done in such a way that the generators of the Rees ideal computed in the { \tt SimisStrategy} phase are not lost; the program computes other generators of the Rees ideal while keeping the generators it found before attaining {\tt HybridLimit}.

There is yet another method for computing the Rees ideal called {\tt SaturationStrategy}. In this option, the whole Rees ideal is computed by saturating the defining ideal of the symmetric algebra with respect to a non-zero element in $R$ (we assume $R$ to be a domain). This strategy appears to be slower in some examples, though one might be able to  improve this option in the future by stopping the computation of the saturation at a certain step.

Computing inverse maps is the most important function of this package and is done by the function {\tt inverseOfMap} (or by running {\tt RationalMap{${}^\wedge$}-1}). According to \autoref{T birationality}, there are two ways to compute the inverse of a map: $(1)$ by finding any syzygy of the Jacobian dual matrix, and $(2)$ by finding a sub-matrix of $\psi$ of rank $\edim(R)-1$. Each way has its benefits. Method $(1)$ is quite fast in many cases, however, method $(2)$ is very useful if the rank of the Jacobian dual matrix  $\psi$ is relatively small compared to the degrees of the entries of $\psi$. Our function {\tt inverseOfMap} starts by using the second method and later switches to the first method if the second method didn't work.  The timing of this transition from the first method to the second method is controlled by the option {\tt MinorsLimit}. Setting {\tt MinorsLimit} to zero will mean that no minors are checked and the inverse map is computed just by looking at the syzygies of $\psi$.  If {\tt MinorsLimit} is left as null (the default value), these functions will determine a value using a heuristic that depends on the varieties involved.

 In addition, to improve the speed of the function {\tt inverseOfMap}, we have two other options, {\tt AssumeDominant} and {\tt CheckBirational}. If {\tt AssumeDominant} is set  to be {\tt true},  then  {\tt inverseOfMap} assumes that the map from $X$ to $Y$ is dominant and does not compute the image of the map; this is time-consuming in certain cases as it computes the kernel of a ring map.  However, this function goes through a call to {\tt idealOfImageOfMap} which first checks whether the ring map is injective (at least if the target is a polynomial ring) using the method described in \cite[Proposition 1.1]{SimisTwoDifferentialThemes}.
 Similarly, if {\tt CheckBirational} sets {\tt false}, {\tt inverseOfMap} will  not check birationality although it still computes the Jacobian dual matrix.  The option {\tt QuickRank} is available to many functions.  At various points, the rank of a matrix is computed, and sometimes it is faster to compute the rank of an interesting-looking submatrix (using the tools of the package {\tt FastMinors}, \cite{FastMinorsArticle}).  Turning {\tt QuickRank} off will make showing that certain maps are birational slower, but will make showing that certain maps are \emph{not} birational faster.  There is a certain amount of randomness in the functions of {\tt FastMinors}, and so occasionally rerunning a slow example will result in a massive speedup.

 In general, as long as {\tt Verbosity} is {\tt >= 1}, the function will make suggestions as to how to run it more quickly.  For example.  
{
{\scriptsize\color{blue}
\begin{verbatim}
    i1 : loadPackage "RationalMaps";

    i2 : Q=QQ[x,y,z,t,u];

    i3 : f = map(Q,Q,matrix{{x^5,y*x^4,z*x^4+y^5,t*x^4+z^5,u*x^4+t^5}});

    o3 : RingMap

    i4 : phi=rationalMapping(f)
                                    5   4    5    4    5    4    5    4
    o4 = Proj Q - - - > Proj Q   {x , x y, y  + x z, z  + x t, t  + x u}

    o4 : RationalMapping

    i5 : time inverseOfMap(phi, CheckBirational=>false, Verbosity => 1);
    inverseOfMapSimis: About to find the image of the map.  If you know the image,
            you may want to use the AssumeDominant option if this is slow.
    inverseOfMapSimis: About to check rank, if this is very slow, you may want to try turning QuickRank=>false.
    inverseOfMapSimis: About to check rank, if this is very slow, you may want to try turning QuickRank=>false.
    inverseOfMapSimis: About to check rank, if this is very slow, you may want to try turning QuickRank=>false.
    inverseOfMapSimis: About to check rank, if this is very slow, you may want to try turning QuickRank=>false.
    inverseOfMapSimis: About to check rank, if this is very slow, you may want to try turning QuickRank=>false.
    inverseOfMapSimis:  We give up.  Using the previous computations, we compute the whole
            Groebner basis of the Rees ideal.  Increase HybridLimit and rerun to avoid this.
    inverseOfMapSimis: Looking for a nonzero minor.
            If this fails, you may increase the attempts with MinorsLimit => #
    inverseOfMapSimis: We found a nonzero minor.
            -- used 0.189563 seconds                                                                     

    o5 : RationalMapping

    i6 : ident = rationalMapping map(Q,Q);

    o6 : RationalMapping

    i7 : o5*phi == ident

    o7 = true
\end{verbatim}
}
Using the {\tt RationalMap{{${}^\wedge$}}-1} syntax to compute inverses of maps will always suppress such output.
{\scriptsize\color{blue}
\begin{verbatim}
    i6 : time phi^-1;
        -- used 0.192791 seconds                                                                     

    o6 : RationalMapping

    i7 : o4 == o7

    o7 = true
\end{verbatim}
}
{\color{black}\normalsize
\section{Embeddings}

Our package also checks whether a rational map $\mathfrak{F} : X \to Y$ is a closed embedding.  The strategy is quite simple.
\begin{enumerate}
\item  We first check whether $\mathfrak{F}$ is regular (by checking if its base locus is empty).
\item  We next invert the map (if possible).
\item  Finally, we check whether the inverse map is also regular.
\end{enumerate}
If all three conditions are met, then the map is a closed embedding and the function returns {\tt true}.  Otherwise, {\tt isEmbedding} returns false.
In the following example which illustrates this, we take a plane quartic, choose a point $Q$ on it, and take the map associated with the divisor $12 Q$.  This map is an embedding by \cite[Chapter IV, Corollary 3.2]{Hartshorne}, which we now verify.
}
{\scriptsize
\color{blue}\begin{verbatim}
    i1 : needsPackage "Divisor"; --used to quickly define a map                                           

    i2 : C = ZZ/101[x,y,z]/(x^4+x^2*y*z+y^4+z^3*x);

    i3 : Q = ideal(y,x+z);

    o3 : Ideal of C

    i4 : f2 = mapToProjectiveSpace(12*divisor(Q));

                            ZZ
    o4 : RingMap C <--- ---[YY ..YY  ]                                                                    
                        101   1    10

    i5 : needsPackage "RationalMaps";

    i6 : time isEmbedding(f2)    
    isEmbedding: About to find the image of the map.  If you know the image,
            you may want to use the AssumeDominant option if this is slow.
    inverseOfMapSimis: About to check rank, if this is very slow, you may want to try turning QuickRank=>false
    inverseOfMapSimis: rank found, we computed enough of the Groebner basis.
            -- used 0.140107 seconds                                                                         

    o6 = true    
\end{verbatim}
}
{\color{black}\normalsize%
\noindent
Notice that {\tt MinorsLimit => 0} by default for {\tt isEmbedding}.  This is because the expressions defining the inverse map obtained from
an appropriate minor frequently are more complicated than the expressions for the inverse map obtained via the syzygies.  Complicated expressions can sometimes slow down the checking of whether the inverse map is regular.

\section{Functionality overlap with other packages}

We note that our package has some overlaps in functionality with other packages.

While the {\tt Parametrization} package \cite{ParametrizationPackage} focuses mostly on curves, it also includes a function called {\tt invertBirationalMap} that has the same functionality as {\tt inverseOfMap}. On the other hand, these two functions were implemented differently and so sometimes one function can be substantially faster than the other.

The package {\tt Cremona} \cite{CremonaSource,CremonaArticle} focuses on very fast probabilistic computation in general cases and very fast deterministic computation for special kinds of maps from projective space. In particular, in {\tt Cremona},

\begin{itemize}
\item{}     {\tt isBirational} gives a probabilistic answer to the question of whether a map between varieties is birational. Furthermore, if the source is projective space, then {\tt degreeOfRationalMap} with {\tt MathMode=>true}  gives a deterministic answer that can be faster than what our package  provides with {\tt isBirationalMap}.
\item{}  {\tt inverseMap} gives a very fast computation of the inverse of a birational map if the source is projective space and the map has maximal linear rank. If this function is passed a map where the domain is not projective space, then it calls a modified, improved version of {\tt invertBirationalMap} originally from {\tt Parametrization}. Even in some cases with maximal linear rank, our {\tt inverseOfMap} function appears to be quite competitive, however.
\end{itemize}

The package {\tt ReesAlgebra} \cite{ReesAlgebraSource, ReesAlgebraArticle} includes a function {\tt jacobianDual} which computes the jacobian dual matrix.  We also have a function {\tt jacobianDualMatrix} which computes a weak form of this same matrix.

\section{Comments and comparisons on function speeds}
\noindent
We begin with a comparison using examples with maximal linear rank where {\tt Cremona} excels.  These examples were executed using version 5.1 of {\tt Cremona} and version 1.0 of {\tt RationalMaps} running Macaulay2 1.19.1.1 on Ubuntu 20.04.

Indeed, in this example (taken from {\tt Cremona}'s documentation), {\tt Cremona} is substantially faster.
{\scriptsize
\color{blue}\begin{verbatim}
    i1 : loadPackage "Cremona"; loadPackage "RationalMaps";

    i3 : ringP20=QQ[t_0..t_20];    

    i4 : phi=map(ringP20,ringP20,{t_10*t_15-t_9*t_16+t_6*t_20,t_10*t_14-t_8*t_16+t_5*t_20,t_9*t_14-t_8*t_15+t_4*t_20,
    t_6*t_14-t_5*t_15+t_4*t_16,t_11*t_13-t_16*t_17+t_15*t_18-t_14*t_19+t_12*t_20,t_3*t_13-t_10*t_17+t_9*t_18-t_8*t_19
    +t_7*t_20,t_10*t_12-t_2*t_13-t_7*t_16-t_6*t_18+t_5*t_19,t_9*t_12-t_1*t_13-t_7*t_15-t_6*t_17+t_4*t_19,t_8*t_12
    -t_0*t_13-t_7*t_14-t_5*t_17+t_4*t_18,t_10*t_11-t_3*t_16+t_2*t_20,t_9*t_11-t_3*t_15+t_1*t_20,t_8*t_11-t_3*t_14
    +t_0*t_20,t_7*t_11-t_3*t_12+t_2*t_17-t_1*t_18+t_0*t_19,t_6*t_11-t_2*t_15+t_1*t_16,t_5*t_11-t_2*t_14+t_0*t_16,
    t_4*t_11-t_1*t_14+t_0*t_15,t_6*t_8-t_5*t_9+t_4*t_10,t_3*t_6-t_2*t_9+t_1*t_10,t_3*t_5-t_2*t_8+t_0*t_10,t_3*t_4
    -t_1*t_8+t_0*t_9,t_2*t_4-t_1*t_5+t_0*t_6});

    o4 : RingMap ringP20 <--- ringP20                                                                     

    i5 : time inverseOfMap(phi, Verbosity=>0);-- Function from "RationalMaps"                              
            -- used 0.118508 seconds                                                                             
    
    o5 : RationalMapping

    i6 : time inverseMap phi;
            -- used 0.0370978 seconds                                                                            
    
    o6 : RingMap ringP20 <--- ringP20                                                                         
    
    i7 : o5 == rationalMapping o6    
    
    o7 = true
\end{verbatim}
}
{\color{black}\normalsize}
However, sometimes the {\tt RationalMaps} function is faster, even in examples with maximal linear rank (a good source of examples where different behaviors can be seen can be found in the documentation of {\tt Cremona}). 
}
{\color{black}\normalsize
We now include an example where the map does not have the maximal linear rank.  
}
{\scriptsize
\color{blue}\begin{verbatim}
    i1 : loadPackage "Cremona"; loadPackage "RationalMaps";

    i3 : Q=QQ[x,y,z,t,u];    

    i4 : phi=map(Q,Q,matrix{{x^5,y*x^4,z*x^4+y^5,t*x^4+z^5,u*x^4+t^5}});    

    o4 : RingMap Q <--- Q                                                                                     

    i5 : (time g = inverseOfMap(phi, Verbosity=>0));
         -- used 0.233111 seconds                                                                             

    i6 : (time f = inverseOfMap(phi, Verbosity=>0, MinorsLimit=>0));
         -- used 60.1969 seconds                                                                          

    i7 : (time h = inverseMap(phi)); -- Function from "Cremona"                                           
         -- used 49.2842 seconds                                                                              

    o7 : RingMap Q <--- Q                                                                                     

    i8 : f == rationalMapping h

    o8 = true

    i9 : g == rationalMapping h

    o9 = true
\end{verbatim}
}
{\color{black}\normalsize
In the previous example, setting {\tt MinorsLimit=>0} makes {\tt inverseOfMap} much slower -- approximately the same speed as the  corresponding command from {\tt Cremona}.
The takeaway for the user should be that changing the options {\tt Strategy}, {\tt HybridLimit}, {\tt MinorSize}, and {\tt QuickRank}, can make a large difference in performance.


We conclude with discussions of the limits of this pacakge.
A work of O.~Gabber shows that if $f : \bP^n \to \bP^n$ is defined by forms of degree $d$, then its inverse can be defined by forms of degree $d^{n-1}$, \cite{BassConnellWrightJacobianConjecture}.  This bound is sharp, as the map
    \[
        (x_0^d : x_1 x_0^{d-1} : x_2 x_0^{d-1} - x_1^d : \dots : x_n x_0^{d-1} - x_{n-1}^d)
    \]
    has inverse given by forms of degree $d^{n-1}$, see \cite{HassanzadehSimisBoundsOnDegreesOfBirat}.  Thus we might expect that this family of maps would be good to explore to see the limits of {\tt RationalMaps}.  We ran these examples with the following code.
}
{\scriptsize
\color{blue}
\begin{verbatim}
    R = ZZ/101[x_0..x_n];
    L = {x_0^d, x_1*x_0^(d-1)} | toList(apply(2..n, i -> (x_i*x_0^(d-1) + x_(i-1)^d)));
    psi = map(R, R, L);
    time inv = inverseOfMap(psi, AssumeDominant=>true, CheckBirational=>false, Verbosity=>0);
\end{verbatim}%
}%
{\color{black}\normalsize
When $n = 3$ (we are working on $\bP^3$) here is a table showing the computation time, in seconds, to find the inverse map for various values of $d$.  The degrees are those we would expect in this example (when $d = 100$, the degree of the forms in the inverse is $10000$).  Note that  {\tt Cremona} has very similar performance for these examples in $\bP^3$ ($n = 3$), but seems substantially slower than {\tt RationalMaps} as we increase the dimension.
\begin{center}
\begin{tabular}{c|c|c|c|c|c|c|c}
  $d$ & 5 & 10 & 20 & 40 & 60 & 80 & 100 \\ \hline
  seconds  &  0.0925 &  0.0958 & 0.1402 &  1.0667 & 7.2652 &  37.4577 & 135.915 \\
\end{tabular}
\end{center}
However, as the size of projective space increases, this becomes much slower.  Here is a table when $n = 4$.
\begin{center}
  \begin{tabular}{c|c|c|c|c|c|c|c|c}
    $d$ & 5 & 8 & 10 & 11 & 12 & 13 & 14 & 15 \\ \hline
    seconds  &  0.1523 & 1.3115 &  7.4682 &  14.9912 &   28.8554 & 57.1229 &   120.778 &  217.706 \\
  \end{tabular}
\end{center}

We conclude with a table when $n = 5$.
\begin{center}
  \begin{tabular}{c|c|c|c|c}
    $d$ & 3 & 4 & 5 & 6  \\ \hline
    seconds  &  0.2619 & 4.8770 &  134.424 &  2713.56  \\
  \end{tabular}
\end{center}
Note the $d = 6$ case took more than 45 minutes.

Finally, Zhuang He and Lei Yang, working under the direction of Ana-Maria Castravet, communicated to us that they used {\tt RationalMaps} to help understand and compute the inverse of a rational map from $\bP^3$ to $\bP^3$, see \cite{HeYangTheMoriDreamSpacePropertyOfBlowups}.  Quoting Zhuang He, this rational map is ``induced by a degree 13 linear system with the  base locus at 6 very general points in P3 and 9 lines through them''.  From a computational perspective, this map was given by 4 degree 13 forms, with 485, 467, 467, and 467 terms respectively.  Computing the inverse of this map took several hours, but it was successful.   

\bibliographystyle{skalpha}
\bibliography{MainBib.bib}
}
\color{black}
\end{document}